\documentclass[11pt]{amsart}
\usepackage{amsfonts}
\usepackage{latexsym}
\usepackage{amssymb}
\usepackage{pdfsync}
\renewcommand{\H}{\mathbb{H}}

\newcommand{\R}{\mathbb{R}}

\newcommand{\cC}{\mathcal{C}}

\newcommand{\cL}{\mathcal{L}}

\newcommand{\cU}{\mathcal U}

\newcommand{\ccI}{\mathfrak{I}}

\renewcommand{\span}{\mbox{\rm span}}

\newcommand{\ep}{\varepsilon}

\newcommand{\ph}{\varphi}

\newcommand{\sm}{\setminus}
\newcommand{\taudot}{\stackrel{.}{\tau}}

\newcommand{\lra}{\longrightarrow}

\newcommand{\der}{\partial}

\newtheorem{The}{Theorem}[section]
\newtheorem{Lem}{Lemma}
\newtheorem{Def}{Definition}

\newtheorem{Cor}{Corollary}[section]

\begin{document}

\title[Intersections of intrinsic submanifolds in the Heisenberg group]
{{\bf Intersections of intrinsic submanifolds in the Heisenberg group}}
\author{Gian Paolo Leonardi}
\address{Gian Paolo Leonardi, Dipartimento di Matematica \\via Campi 213/b\\ I-41100, Modena}
\email{gianpaolo.leonardi@unimore.it}
\author{Valentino Magnani}
\address{Valentino Magnani, Dipartimento di Matematica \\
Largo Pontecorvo 5 \\ I-56127, Pisa}
\email{magnani@dm.unipi.it}
\begin{abstract}
In the first Heisenberg group, we show that the intersection of two intrinsic 
submanifolds with linearly independent horizontal normals
locally coincides with the image of an injective continuous curve. 
The key tool is a chain rule that relies on a recent result by Dafermos.
\end{abstract}

\maketitle

\tableofcontents{}

\bigskip
\footnoterule{
The authors have been supported by GNAMPA of INdAM:
project ``Maps, measures and nonlinear calculus in
stratified groups, 2008'' 
{\em keywords:} Heisenberg group, intrinsic regular surfaces,
implicit function theorem}

\pagebreak

\section{Introduction}

This work is mainly concerned with low dimensional intrinsic regular
sets in the first Heisenberg group $\H$, equipped with a sub-Riemannian distance. 
Some motivations for this study stem from the project of developing 
Geometric Measure Theory in stratified groups.

The notion of ``regular set'' here refers to the sub-Riemannian metric structure.
In the family of Heisenberg groups $\H^n$, a natural class of intrinsic 
regular surfaces  has been found
by B. Franchi, R. Serapioni and F. Serra Cassano, \cite{FSSC6}.
The choice of these sets naturally arises from the so-called ``Rumin complex'',
\cite{Rum90}, according to which $\H$-regular surfaces of \cite{FSSC6}
can be seen as ``regular currents'' defined on the space of compactly 
supported Rumin differential forms. This complex singles out 
two distinct classes of $\H$-regular surfaces. 
The {\em low dimensional $\H$-regular surfaces} are $C^1$ smooth horizontal
submanifolds, whose topological dimension is less than or equal to $n$. 
On the opposite side, the {\em low codimensional $\H$-regular surfaces}
have topological dimension greater than $n$ and are not differentiable in the classical sense.
Recently, G. Arena and R. Serapioni have characterized all
these sets as {\em intrinsic differentiable} intrinsic graphs,
\cite{AreSer09}.

A low codimensional $\H$-regular surface is locally defined
as a level set of an $\R^k$-valued differentiable mapping
on an open set of $\H^n$, where the differential is surjective and
differentiability is meant with respect to intrinsic dilations
and the group operation, see Section~\ref{basicf}. 
It has been proved in \cite{FSSC6} that this set, 
although it may not be $C^1$ smooth in the classical sense,
is locally an intrinsic graph and its
$(2n$$+$$2$$-$$k)$-dimensional Hausdorff measure
with respect to the sub-Riemannian metric structure
can be computed by an area-type formula.
It is worth to stress that these results are achieved when $k\leq n$, 
since the corresponding semidirect factorizations of $\H^n$ can be used
to realize the intrinsic graph structure of the level set.
A consequence of the above mentioned area-type formula 
is the sub-Riemannian coarea formula for Lipschitz mappings 
$f:A\subset\H^n\lra\R^k$, under the condition $k\leq n$, \cite{Mag15}.
The validity of this formula for $k>n$ is still not clear.
As suggested by this problem, one is lead to consider level sets
of $\R^k$-valued differentiable mappings with surjective differential
in the case $k>n$. These level sets do not belong to
the class of $\H$-regular surfaces of \cite{FSSC6}.

Our aim is to understand the structure of these sets in the simplest situation, 
namely, in the first Heisenberg group $\H$. This amounts to study whether 
a possible ``implicit curve theorem'' for $\R^2$-valued differentiable mappings
is available. An equivalent formulation of this fact
is to study the intersection of two $\H$-regular surfaces of $\H$
having linearly independent horizontal normals at the intersection points.
Our main result is the following
\begin{The}\label{ImplCurv}
Let $\Sigma_1$ and $\Sigma_2$ be two $\H$-regular surfaces 
of $\H$, where $p\in\Sigma_1\cap\Sigma_2$ and the horizontal normals
${\nu_1}_H(p)$ and ${\nu_2}_H(p)$ are linearly independent. 
Then there exists a neighbourhood $O\subset\H$ of $p$ such that the set
$\Sigma_1\cap\Sigma_2 \cap O$ coincides with the image of an 
injective continuous curve $\Gamma:[0,1]\to\H$. 
\end{The}
The proof of this theorem differs from both the Euclidean case
and the Heisenberg case considered in \cite{FSSC6}.
In the sequel, we present the main ideas that lead to this result.
It is clearly not restrictive to assume that $p$ is the origin and that 
both $\Sigma_1$ and $\Sigma_2$ are defined around this point as zero level 
sets of differentiable functions $f_1$ and $f_2$, respectively.
Notice that in particular $\Sigma_2$ need not be a graph in the classical sense,
but it can be injectively parametrized by an intrinsic graph mapping
$\Phi_2:n\to n\cdot\phi(n)$ that involves the group operation, where
$n$ and $\phi(n)$ belong to a vertical subgroup and to a horizontal subgroup of $\H$,
respectively. This was observed in \cite{FSSC3} and subsequently developed 
for higher codimensions in \cite{FSSC6}.
Then we notice that the intersection of $\Sigma_1$ and $\Sigma_2$ around the
origin amounts to the $\Phi_2$-image of the zero level set of $f_1\circ\Phi_2$.
Both $f_1$ and $\Phi_2$ are not differentiable in the classical sense, hence
the level set of their composition may have a priori a highly nontrivial structure.
The classical Dini's scheme to prove the implicit function theorem
suggests us to find directions, or more precisely a family of curves,
along which $f_1\circ\Phi_2$ is strictly monotone.
However, in our case these curves are solutions to a continuous ODE, hence
the corresponding flow is not uniquely defined. By Lemma~\ref{lemMajer},
that has been kindly pointed out to us by P. Majer, 
we select a continuous and ``monotone'' flow
of solutions that leads us to Theorem~\ref{teo:curva}.
Here we show that whenever this family of curves is available, 
then there exists an injective continuous parametrization of the level set.
To find these curves, we follow the recent theory developed in \cite{ASCV},
\cite{BSC1} and \cite{BSC2}.
First of all, the work of L. Ambrosio, F. Serra Cassano and D. Vittone,
\cite{ASCV}, shows that $\ph$, as a a scalar function, is a distributional solution of the Burgers equation
\begin{equation}\label{Burgers}
\frac{\der}{\der y} \ph+\det C\,\frac{\der}{\der t}\big(\ph^2\big)=
-\frac{(Yf_2)\circ\Phi_2}{(Xf_2)\circ\Phi_2}\,.
\end{equation}
Here we are following notation and terminology of both Section~\ref{basicf}
and Section~\ref{DirectDerivChRu}.
In the recent work \cite{BSC2}, F. Bigolin and F. Serra Cassano
establish the full characterization of $\H$-regular surfaces as intrinsic graphs
of distributional solutions of the Burgers equation. Here the authors find 
an interesting connection with a recent result of C. Dafermos, \cite{Dafermos06}.

In short, all of these results show that the required family of curves 
exactly corresponds to the characteristics of \eqref{Burgers}
in the vertical plane. In fact, these characteristics are precisely 
those curves that are ``lifted'' by $\Phi_2$ to $C^1$ smooth horizontal 
curves in $\Sigma_2$. 
This is somewhat surprising, since $\Sigma_2$ is not even Lipschitz regular from the Euclidean viewpoint.
Then the claim of Theorem~\ref{ImplCurv} is established if we show that
$f_1\circ\Phi_2$ is differentiable along these curves and it has 
nonvanishing derivative. 

In Theorem~\ref{Taytyp}, taking into account that differentiability of $f_1$
implies classical differentiability along horizontal directions and using a
recent result by G. Arena and R. Serapioni, \cite{AreSer09},
we have found an {\em intrinsic Taylor expansion} of $f_1\circ\Phi_2$.
This theorem implies in turn a partial differentiability along directions of 
characteristics, see Corollary~\ref{directderiv}.
Nevertheless, this is still not enough to obtain differentiability
of $f_1\circ\Phi_2$ along characteristics, since $f_1\circ\Phi_2$ is not
regular in the Euclidean sense. The final step to overcome this problem
is to join Theorem~\ref{Taytyp} with the above mentioned result of C. Dafermos,
in the version stated in \cite{BSC2}.

As a consequence, in Theorem~\ref{chainr} we establish a chain rule for
the composition $f_1\circ\Phi_2$, where  
the linear independence of $\nabla_Hf_1(0)$ and $\nabla_Hf_2(0)$ finally
shows that the derivative of $f_1\circ\Phi_2$ along characteristics is
also nonvanishing, leading us to the required strict monotonicity.
Finally, in Subsection~\ref{intLip} we add an observation about the 
regularity of the implicit curve. We show that the
implicit curve, as a set, has a sort of cone-type Lipschitz continuity, in 
analogy with the intrinsic Lipschitz continuity introduced in \cite{FSSC7}. 
Here the main difference is that the intrinsic cone refers to a 
factorization of the Heisenberg group that is no longer a semidirect product.

\vskip.1cm
\noindent
{\bf Acknowledgements.}
We gratefully thank Pietro Majer for suggesting us a useful lemma. 
We are also indebted with Francesco Serra Cassano for fruitful discussions
and for having informed us on his recent results with Francesco Bigolin.

\section{Basic facts}\label{basicf}

We represent the Heisenberg group $\H$ as a 3-dimensional Hilbert space,
equipped with orthogonal subspaces $H_1$ and $H_2$ such that $\H=H_1\oplus H_2$,
$\dim(H_1)=2$ and $\dim(H_2)=1$. 
Any element $x\in\H$ can be uniquely represented by $x=x_1+x_2$, with $x_j\in H_j$
and $j=1,2$. Let us now fix an orthonormal basis $(e_1,e_2,e_3)$
of $\H$, where $e_3\in H_2$. Taking into
account the above direct sum, we will also
use the further decomposition $x_1 = x_{1,1} e_1+ x_{1,2} e_2\in H_1$.  
The group operation in $\H$ is defined as follows: for any $x,y\in \H$
\begin{equation}\label{group}
x\cdot y=x+y+\omega(x_1,y_1)\,,
\end{equation}
where $\omega:H_1\times H_1\lra H_2$, $\omega(x_1,y_1)=\bar\omega(x_1,y_1)\, e_3$ and 
\[
\bar\omega(x_1,y_1)= x_{1,1}\,y_{1,2}-y_{1,1}\,x_{1,2}.
\]
The so-called ``intrinsic dilations'' of $\H$ are defined as $\delta_r(x)=rx_1+r^2x_2$.
We fix the metric structure on $\H$ introducing the
{\em homogeneous norm} 
\[
\|x\|=\max\left\{|x_1|,\sqrt{|x_2|}\right\}\,,
\] 
that satisfies the triangle inequality
$\|x\cdot y\|\leq\|x\|+\|y\|$ and the homogeneity
$\|\delta_rx\|=r\|x\|$ for all $x,y\in\H$ and $r>0$. This homogeneous norm
defines a distance on $\H$, setting $d(x,y) = \|x^{-1}\cdot y\|$ for all $x,y\in \H$. 
With respect to this distance, the closed ball of center $x$ and radius $r>0$ is denoted by
$D_{x,r}$. The point $x$ is omitted when it coincides with the origin.
If $E\subset\H$, we also use the notation $D^E_{x,r}=D_{x,r}\cap E.$
We use the same conventions for the open balls $B_{x,r}$ of center $x$ and radius $r$.
In sum, we wish to stress that the previous assumptions on $\H$ allow us
to regard it as a Hilbert space, a Lie group with group operation \eqref{group} and
a metric space equipped with distance $d$, simultaneously.

Throughout the paper $\Omega$ will denote an open subset of $\H$.
A mapping $f:\Omega\lra\R^2$ is {\em differentiable} at $x\in\Omega$
if the following Taylor expansion 
\[
f(y)=f(x)+\nabla_Hf(x)(x^{-1}y)+o\big(d(x,y)\big)
\]
holds as $d(x,y)\to0$. From this definition, it follows that
$\nabla_Hf(x):\H\lra\R^2$ is also {\em homogeneous}, namely, 
$\nabla_Hf(x)(\delta_rz)=r\nabla_Hf(x)(z)$ for all $z\in\H$ and $r>0$.
Let $f_1$ and $f_2$ denote the components of $f$. Then we will also think of
$\nabla_Hf(x)$ as a $2\times 3$ matrix of rows 
$\nabla_Hf_1(x)$ and $\nabla_Hf_2(x)$. 
These vectors are the well known {\em horizontal gradients} of the components
$f_1$ and $f_2$, respectively.
Notice that the left invariant vector fields 
$(X_1,X_2)$ spanned by the orthonormal vectors
$(e_1,e_2)$ yield the form of $\nabla_Hf(x)$ 
as a matrix with respect to the fixed scalar product, namely,
\[
\big(\nabla_Hf(x)\big)^i_j=X_jf_i(x) \text{ and } \big(\nabla_Hf(x)\big)^i_{3}=0 
\]
for $i,j=1,2$. 
We denote by $\cC^1(\Omega,\R^2)$ the linear space of mappings 
$f:\Omega\lra\R^2$ that are differentiable and such that
$x\lra\nabla_Hf(x)$ is continuous.
Taking into account the ``stratified mean value theorem'' (1.41) of \cite{FS82},
the next lemma can be shown by standard arguments.
\begin{Lem}\label{unidiff}
Let $f\in\cC^1(\Omega,\R^2)$. Then for every $p\in\Omega$ there exists $\delta>0$ 
and a nonincreasing function $\omega_p:[0,\delta[\lra\R$ infinitesimal
at zero and depending on the modulus of continuity of $x\lra\nabla_Hf(x)$,
such that $D_{p,\delta}\subset\Omega$ and for all $0<t\leq\delta$, we have
\[
\sup\left\{\frac{\left|f(y)-f(x)-\nabla_Hf(x)(x^{-1}\cdot y)\right|}
{d(x,y)}:x,y\in D_{p,\delta},\; 0<d(x,y)\leq t\right\}\leq\omega_p(t)\,.
\]
\end{Lem}
\begin{Def}{\rm
We say that any one dimensional subspace $H$ of $\H$, that is contained in $H_1$,
is a {\em horizontal subgroup}. We say that $N$ is a {\em vertical subgroup} if
$N=H+H_2$, where $H$ is a horizontal subgroup $H$.
}\end{Def}
\begin{Def}\rm
We say that $\Sigma\subset\H$ is an {\em $\H$-regular surface}
if for each point $p\in\Sigma$ there exists an open set $\cU\subset\H$
containing $p$ and a function $f\in\cC^1(\Omega,\R)$ such that
$\Sigma\cap\cU=f^{-1}(0)$ and $\nabla_H f(q)\neq0$ for all $q\in\cU$.
A {\em horizontal normal} of $\Sigma$ at $p$ is given by $\displaystyle 
\nu_H(p)=\frac{\nabla_Hf(p)}{|\nabla_Hf(p)|}$.
\end{Def}
\begin{Def}\label{Jsemid}\rm
Let $H$ be a horizontal subgroup and let $N$ be a normal subgroup such that $H\oplus N=\H$.
Then we define the {\em factorizing mapping} $J$ associated to this
direct sum by $J(n,v)=n\cdot v$, where $J:N\times H\lra\H$.
A direct computation shows that the factorizing mapping $J$ is a diffeomorphism 
and its inverse defines the {\em canonical projections} 
$\pi_N:\H\lra N$, $\pi_H:\H\lra H$ by the formula $J^{-1}(x)=\big(\pi_N(x),\pi_H(x)\big)$.
Explicitly, we have
\[
\pi_N(x)=x-x_H-\omega(x_1,x_H)  \quad\mbox{and}\quad \pi_H(x)=x_H,
\]
where $x_H\in H$ and $x-x_H\in N$.
\end{Def}
We will also use some recent results from the theory of scalar conservation laws.
\begin{Def}\rm
Let $O$ be an open set of $\R^2$, let $\beta\neq0$ and let $\ph:O\lra\R$
be continuous.
We consider a bounded measurable function 
$g:I\times J\lra\R$, where $I$ and $J$ are open intervals
and $I\times J\subset O$. Let $[a,b]\subset I$ and let $\tau\in C^1([a,b],J)$.
We say that $\tau$ is a {\em characteristic} associated to the Burgers equation
\begin{equation}\label{CcharDaf}
\frac{\der}{\der y} \ph+\frac{\beta}{2}\,\frac{\der}{\der t}\big(\ph^2\big)=g
\end{equation}
if $\dot\tau(y)=\beta\,\ph\big(y,\tau(y)\big)$ for all $y\in(a,b)$.
\end{Def}

\medskip
\noindent
Following Theorem~2.1 of \cite{BSC2}, the next theorem states
a version of a result by C. M. Dafermos, \cite{Dafermos06},
that is one of the key points to establish the chain rule of Theorem~\ref{chainr}.
\begin{The}\label{Dafermos}
Let $I$ and $J$ be two open intervals and let 
$g:I\times J\lra\R$ be bounded measurable such that $g(\eta,\cdot)$
is continuous on $J$ for any $\eta\in I$.
Let $[a,b]\subset I$ and let $\tau\in C^1\big([a,b],J)$
be a characteristic associated with a continuous distributional solution
$\ph$ of \eqref{CcharDaf} with $\beta=1$.
We set $\nu(\eta)=\ph\big(\eta,\tau(\eta)\big)$ for any $\eta\in[a,b]$.
Then $(\tau,\nu)$ satisfies the system
\begin{eqnarray}
\left\{\begin{array}{l}
\dot\tau(s)=\nu(s) \\
\dot\nu(s)=g(s,\tau(s)) 
\end{array}\right.\quad \mbox{for all}\quad s\in[a,b]\,.
\end{eqnarray}
In particular, $\dot\tau$ is Lipschitz continuous on $[a,b]$.
\end{The}
%
%
%
%
%
%
%
%
%
%
%
%
%
\section{Taylor-type expansion and chain rule}\label{DirectDerivChRu}
%
%
%
%
%

In order to simplify the statements of the main results of this section, 
we first fix a number of basic assumptions we will use throughout. 
First, a vertical subgroup $N\subset\H$ along with
a horizontal subgroup $H\subset\H$ such that $H\oplus N=\H$ will be fixed.
 
We will consider a unit (horizontal) vector $b_1\in H$ and an orthonormal basis 
$(b_2,e_3)$ of $N$, where $b_1$ and $b_2$ span $H_1$, although they need not be orthogonal.
We make explicit the change of variable with respect to the orthonormal basis $(e_1,e_2)$ of $H_1$,
setting 
\[
b_1=c_1^1\,e_1+c_1^2\,e_2,\quad b_2=c_2^1\,e_1+c_2^2\,e_2 \quad \mbox{and}\quad 
C=\left(\begin{array}{cc} c^1_1 & c^1_2 \\ c^2_1 & c^2_2 \end{array}\right).
\]
Then the corresponding left invariant vector fields associated to $b_1$ and $b_2$ are 
\[
Y_1=c_1^1\,X_1+c_1^2\,X_2\quad \mbox{ and }\quad Y_2=c_2^1\,X_1+c_2^2\,X_2.
\]
We fix a function $f\in\cC^1(\Omega,\R^2)$ of components $f_1$ and $f_2$, 
a point $x_0\in\Omega$ where $\nabla_Hf(x_0)$ is surjective and assume
\begin{equation}\label{surhdi}
Y_1f_2(x_0)\neq0.
\end{equation}
This is not restrictive, up to exchanging the components of $f$.
The point $x_0$ can be written in a unique way by the group product 
$x_0=n_0\cdot v_0$, where $n_0\in N$ and $v_0\in H$.
Taking into account Definition~\ref{Jsemid}, we can find $s>0$ such that 
\[
B^N_{n_0,s}\cdot B^H_{v_0,s}:=J\big(B^N_{n_0,s}\times B^H_{v_0,s}\big)\subset\Omega,
\]
where $J$ is the factorizing mapping associated to the direct sum $H\oplus N=\H$.
It is also not restrictive to assume that $Y_1f_2\neq0$ everywhere on $B^N_{n_0,s}$.

Condition \eqref{surhdi} suffices to apply the implicit function theorem of \cite{FSSC3},
getting an ``intrinsic graph mapping'' $\Phi_2:B^N_{n_0,s}\lra \H$, for some possibly 
smaller $s>0$, such that 
\[
\Phi_2(n)=n\cdot\phi_2(n)\quad\mbox{ and }\quad \phi_2:B^N_{n_0,s}\lra H
\]
and where the following conditions
\[
\Phi_2(n_0)=x_0\quad\mbox{ and } \quad f_2\big(\Phi_2(n)\big)=f_2(x_0)
\]
are satisfied for all $n\in B^N_{n_0,s}$. Finally, we introduce the uniquely defined mapping
\[
\ph_2:B^N_{n_0,s}\lra\R, \quad \mbox{such that}\quad \phi_2(n)=\ph_2(n)\,b_1\,.
\]

In the sequel, the previous assumptions will be understood.
\begin{The}[Taylor-type expansion]\label{Taytyp}
Let $\bar n=\bar \eta_2b_2+\bar \tau_1e_3\in B^N_{n_0,s}$ be fixed and consider
$\bar v=\phi_2(\bar n)\in H$, along with
$
\bar x:=\Phi_2(\bar n)=\bar n\cdot \bar v=\bar \eta_1b_1+\bar \eta_2\,b_2+\bar \tau e_3\in\H.
$
Then the following Taylor-type expansion holds
\begin{equation}\label{Taylorintr}
f_1\circ\Phi_2(n)=f_1\circ\Phi_2(\bar n)-\frac{(\eta-\bar \eta_2)}{Y_1f_2(\bar x)}\;\det
\left(\begin{array}{cc}
Y_1f_1(\bar x) & Y_2f_1(\bar x)
\\ Y_1f_2(\bar x) & Y_2f_2(\bar x)
\end{array}\right)+o(\|(\eta-\bar \eta_2)\,b_2+\tau'\,e_3\|)\,,
\end{equation}
where $n=\eta\, b_2+\tau\, e_3$ and $\tau'$ is given by the change of variable 
\begin{equation}\label{changevt}
\tau'=\tau-\bar\tau-2\eta\,\bar \eta_1\,\det(C)+
\bar \eta_1\,\bar \eta_2\,\det C\,.
\end{equation}
\end{The}
{\sc Proof.}
We set $g_1(\xi)=f_1(\bar x\cdot\xi)$ and observe that 
\begin{eqnarray*}
f_1\big(n\cdot\phi_2(n)\big)=g_1\big(\bar x^{-1}\cdot n\cdot\phi_2(n)\big)
=g_1\big(u\cdot (\phi_2)_{\bar x^{-1}}(u)\big)\,,
\end{eqnarray*}
where $(\phi_2)_{\bar x^{-1}}$ is the ``translated function'', 
as introduced in (i) of Proposition~3.6 in \cite{AreSer09}.
Then $u$ is given by the expression
\[
\tau_{\bar x^{-1}}(n):=\bar x^{-1}\cdot n\cdot\bar x\cdot\pi_N(\bar x^{-1})\,,
\]
where the projection $\pi_N$ is recalled in Definition~\ref{Jsemid}.
Notice that $\tau_{\bar x^{-1}}(N)\subset N$ and 
\begin{eqnarray*}
\tau_{\bar x^{-1}}(n)&=&
\bar x^{-1}\cdot n\cdot\bar x\cdot\pi_N(\bar x^{-1})\cdot \pi_H(\bar x^{-1})\cdot 
\big(\pi_H(\bar x^{-1})\big)^{-1} \\
&=& \bar x^{-1}\cdot n\cdot \big(\pi_H(\bar x^{-1})\big)^{-1}\,.
\end{eqnarray*}
Taking into account the expression of $\pi_H$ in Definition~\ref{Jsemid}
and that $\bar x^{-1}=-\bar x$, we get 
\[
\big(\pi_H(\bar x^{-1})\big)^{-1}=(\bar x)_H=\pi_H(\bar x)=\bar v\,.
\]
Then the following formulae hold
\[
\tau_{\bar x^{-1}}(n)=\bar x^{-1}\cdot n\cdot \bar v=
c_{\bar v^{-1}}(\bar n^{-1}n),
\]
where $c_{\bar v^{-1}}(x)=\bar v^{-1}\cdot x\cdot \bar v$ is a group isomorphism of $N$. 
Then $\tau_{\bar x^{-1}}(B^N_{n_0,s})=c_{\bar v^{-1}}\big(B^N_{\bar n^{-1}\cdot n_0,s}\big)$.
Notice that this open set contains the origin, hence the new variable 
$u\in c_{\bar v^{-1}}\big(B^N_{\bar n^{-1}\cdot n_0,s}\big)$
varies in an open neighbourhood of the origin.

Let us set the ``translated'' variables $(\eta',\tau')$, such that $u=\eta'b_2+\tau'e_3$.
Taking into account that $u=\tau_{\bar x^{-1}}(n)$,
a direct computation shows that
\begin{equation*}
\left\{\begin{array}{l}
\eta'=\eta-\bar \eta_2 \\
\tau'=\tau-\bar \tau-2\eta\,\bar \eta_1\,\det(C)+
\bar \eta_1\,\bar \eta_2\,\det C
\end{array}\right.\,.
\end{equation*}
Notice that the second equation yields \eqref{changevt}. Now,
we wish to study the local expansion of $f_1\circ\Phi_2$ around zero, 
with respect to the new variables $(\eta',\tau')$. First of all, 
condition \eqref{surhdi} implies that the level set $f_2$ is a low 
codimensional $\H$-regular surface, that is given by its intrinsic graph 
mapping $\Phi_2:B^N_{n_0,s}\lra \H$, with $\Phi_2(n)=n\cdot\phi_2(n)$.
Thus, Theorem~4.2 of \cite{AreSer09} implies that $\phi_2$ is (uniformly)
intrinsic differentiable. According to Definition~3.13 of \cite{AreSer09}, 
intrinsic differentiability of $\phi_2$ corresponds to differentiability 
at the origin of its translated functions, with respect to intrinsic linear 
mappings. In particular, our translated function $(\phi_2)_{\bar x^{-1}}$ satisfies
\begin{equation}\label{intr}
\lim_{u\to0}\frac{\|L_2(u)^{-1}\cdot(\phi_2)_{\bar x^{-1}}(u)\|}{\|u\|}=0\,,
\end{equation}
where $L_2:N\lra H$ is an intrinsic linear mapping. 
Recall from Proposition~3.23(ii) of \cite{AreSer09} that any intrinsic linear mapping
is H-linear. Then $L_2$ is H-linear, i.e.  
a group homomorphism, satisfying $L(\delta_ru)=rL(u)$ for all $u\in N$ and $r>0$.
As a consequence, $L_2$ is a linear mapping satisfying $L_2(\eta' b_2+\tau'e_3)=\alpha\,\eta'b_1$ for 
some fixed $\alpha\in\R$.

We define $g_2(\xi)=f_2(\bar x\cdot \xi)$. Due to the differentiability of $g_2$ at $0$ and
taking into account \eqref{intr}, the chain rule gives
\[
\nabla_Hg_2(0)(u\cdot L_2(u))=\nabla_Hf_2(x_0)(u\cdot L_2(u))=0\,.
\] 
Then the previous equation gives
\[
L_2(u)=-\eta'\,\frac{Y_2f_2(\bar x)}{Y_1f_2(\bar x)}\,b_1\,.
\]
Condition \eqref{intr} implies that $(\phi_2)_{\bar x^{-1}}(u)=L_2(u)\cdot\ep_2(u)$,
where $\|\ep_2(u)\|/\|u\|\to 0$ as $u\to0$.
Setting $L_1=\nabla_Hf_1(\bar x)$, differentiability of $g_1$ at zero yields
\[
g_1(x)-g_1(0)=L_1(x)+\ep_1(x)
\]
where $\|\ep_1(x)\|/\|x\|\to0$ as $x\to0$. It follows that
\begin{eqnarray}\label{IntrExp}
g_1(u\cdot(\phi_2)_{\bar x^{-1}}(u))-g_1(0)&=&L_1(u\cdot(\phi_2)_{\bar x^{-1}}(u))+
\ep_1(u\cdot(\phi_2)_{\bar x^{-1}}(u))\\
&=&L_1(u)+(L_1\circ L_2)(u)+L_1(\ep_2(u))+\ep_1(u\cdot(\phi_2)_{\bar x^{-1}}(u))\nonumber\\
&=&L_1(u)+(L_1\circ L_2)(u)+o(u)\nonumber\\
&=&\eta'\,Y_2f_1(\bar x)-\eta' \frac{Y_2f_2(\bar x)}{Y_1f_2(\bar x)}\,Y_1f_1(\bar x)
+o(\|\eta'\,w+\tau'\,e_3\|)\nonumber
\end{eqnarray}
where $\|o(u)\|/\|u\|\to0$ as $u\to0$, since $\|\ep_1(u\cdot(\phi_2)_{\bar x^{-1}}(u))\|/\|u\|\to0$ as $u\to0$.
This proves the Taylor expansion \eqref{Taylorintr}. $\Box$
\begin{Cor}[Directional derivatives]\label{directderiv}
Under the assumptions of Theorem~\ref{Taytyp},
$f_1\circ\Phi_2$ is partially differentiable at $\bar n$ along
$
\bar z=b_2+2\,\bar \eta_1\,(\det C)\,e_3\in N
$
and there holds
\begin{equation}\label{derivativeform}
\der_{\bar z}\big(f_1\circ\Phi_2\big)(\bar n)
=-\frac{1}{Y_1f_2(\bar x)}\;\det
\left(\begin{array}{cc}
Y_1f_1(\bar x) & Y_2f_1(\bar x) \\ 
Y_1f_2(\bar x) & Y_2f_2(\bar x)
\end{array}\right)\,.
\end{equation}
\end{Cor}
{\sc Proof.} 
The idea is of restricting the expansion \eqref{Taylorintr} to the set of 
points $n=\eta b_2+\tau e_3$ such that the corresponding change of variable
$\eta'b_2+\tau'e_3=\tau_{\bar x^{-1}}(\eta b_2+\tau e_3)$ satisfies $\tau'=0$.
We wish to read this constraint with respect to the initial variables 
$(\eta,\tau)$, where $n=\eta b_2+\tau e_3$.
Recall that the change of variable $u=\eta' b_2+\tau'e_3=\tau_{\bar x^{-1}}(n)$ gives
\begin{equation*}
\left\{\begin{array}{l}
\eta'=\eta-\bar \eta_2 \\
\tau'=\tau-\bar \tau-2\eta\,\bar \eta_1\,\det(C)+
\bar \eta_1\,\bar \eta_2\,\det C
\end{array}\right.\,.
\end{equation*}
Then our constraint $\tau'=0$, yields 
\begin{eqnarray*}
\tau&=&\bar \tau+2\eta\,\bar \eta_1\,\det(C)-\bar \eta_1\,\bar \eta_2\,\det C \\
&=& \bar \tau+2(\eta-\bar\eta_2)\,\bar \eta_1\,\det(C)+\bar \eta_1\,\bar \eta_2\,\det C \\
&=& \bar \tau+2\eta'\,\bar \eta_1\,\det(C)+\bar \eta_1\,\bar \eta_2\,\det C\,.
\end{eqnarray*}
Then we get a line $l(\eta')$ in $N$ of coordinates
\begin{equation}
\left\{\begin{array}{l}
l_1(\eta')=\eta'+\bar \eta_2 \\
l_2(\eta')=\bar\tau+2\eta'\,\bar \eta_1\,\det(C)+\bar \eta_1\,\bar \eta_2\,\det C
\end{array}\right.\,.
\end{equation}
The equation $\bar n\cdot \bar v=\bar x$ yields 
$
\bar n=\bar \eta_2b_2+\bar \tau_1e_3=\bar \eta_2b_2+(\bar \tau+\bar \eta_1\bar \eta_2\det C)\,e_3\,,
$
then 
\[
l_2(\eta')=\bar\tau_1+2\eta'\bar\eta_1\det C\,.
\]
It follows that $l$ can be written as follows
\[
l(\eta')=(\bar\eta_2+\eta')b_2+(\bar\tau_1+2\eta'\bar\eta_1\det C)e_3,
\]
where $l(0)=\bar n$. Thus, due to expansion \eqref{Taylorintr}, we can establish
\begin{equation*}
f_1\circ\Phi_2\big(l(\eta')\big)-f_1\circ\Phi_2\big(l(0)\big)=
-\frac{\eta'}{Y_1f_2(\bar x)}\;\det
\left(\begin{array}{cc}
Y_1f_1(\bar x) & Y_2f_1(\bar x) \\ 
Y_1f_2(\bar x) & Y_2f_2(\bar x)
\end{array}\right)
+ o(\eta')\,,
\end{equation*}
that leads us to the conclusion. $\Box$

\medskip

In the next theorem, we identify $N$ with $\R^2$, through the isomorphism 
$(\eta,t)\lra\eta\, b_2+t\,e_3$.
\begin{The}[Chain rule]\label{chainr}
Let $I$ and $J$ be two open intervals such that $I\times J\subset B^N_{n_0,s}$ 
and define $[a,b]\subset I$.
Let $\tau\in C^1([a,b],J)$ be any characteristic of 
\begin{equation}\label{CcharDafappl}
\frac{\der}{\der \eta}\, \ph_2+\det(C)\,\frac{\der}{\der t}\big(\ph_2^2\big)
=-\frac{Y_2f_2\big(\Phi_2(\eta\, b_2+te_3)\big)}{Y_1f_2(\Phi_2(\eta\, b_2+te_3))}.
\end{equation}
If $\gamma(\eta)=\eta\,b_2+\tau(\eta)\,e_3$, then the composition
$f_1\circ\Phi_2\circ\gamma$ is everywhere differentiable and
\begin{equation}\label{chainrule}
\frac{d}{d\eta}\big(f_1\circ\Phi_2\circ\gamma\big)(\eta)=
-\frac{1}{Y_1f_2\big(\Phi_2\circ \gamma\big)(\eta)}\;\det
\left(\begin{array}{cc}
Y_1f_1\big(\Phi_2\circ \gamma(\eta)\big) & Y_2f_1\big(\Phi_2\circ \gamma(\eta)\big)\\ 
Y_1f_2\big(\Phi_2\circ\gamma(\eta)\big) & Y_2f_2\big(\Phi_2\circ \gamma(\eta)\big)
\end{array}\right)\,.
\end{equation}
\end{The}
{\sc Proof.}
We have to show the differentiability of $f_1\circ\Phi_2\circ\gamma$ 
and the everywhere validity of \eqref{chainrule}. 
Following the notation of Theorem~\ref{Taytyp}, we fix $\bar \eta_2\in[a,b]$, along with
\[
\gamma(\bar \eta_2)=\bar n=\bar \eta_2b_2+\bar \tau_1e_3\quad\mbox{and}\quad 
\tau(\bar \eta_2)=\bar \tau_1.
\]
Taking into account the constraint $\bar x=\bar n\cdot\bar v$, we get
\[
\bar x=(\bar\eta_2b_2+\bar\tau_1e_3)\cdot(\bar\eta_1b_1)=\bar\eta_1b_1+\bar\eta_2b_2+
(\bar\tau_1-\bar\eta_1\bar\eta_2\det C)e_3
\]
where $\bar\eta_1=\ph_2(\bar n)$. It follows that
\begin{equation}\label{tau1}
\bar\tau=\bar\tau_1-\bar\eta_1\bar\eta_2\det C.
\end{equation}
In view of expansion \eqref{Taylorintr}, we get
\begin{equation*}
f_1\circ\Phi_2\big(\gamma(\eta)\big)-f_1\circ\Phi_2(\bar n_2)=
-\frac{(\eta-\bar \eta_2)}{Y_1f_2(\bar x)}\;\det
\left(\begin{array}{cc}
Y_1f_1(\bar x) & Y_2f_1(\bar x)\\ 
Y_1f_2(\bar x) & Y_2f_2(\bar x)
\end{array}\right) 
+o(\|(\eta-\bar \eta_2)\,w+\tilde \tau(\eta)\,e_3\|)\,,
\end{equation*}
where the change of variable \eqref{changevt} gives
$\tilde\tau(\eta)=\tau(\eta)-\bar\tau-2\eta\bar\eta_1\det C+\bar\eta_1\bar\eta_2\det C$, then
\begin{equation}
\tilde\tau(\eta)=\tau(\eta)-\bar \tau_1-2(\eta-\bar \eta_2)
\ph_2(\gamma(\bar \eta_2))\,\det C\,.
\end{equation}
as a consequence of \eqref{tau1}.
Hence differentiability follows if we show that $|\tilde\tau(\eta)|=O(|\eta-\bar \eta_2|^2)$.
Then we consider
\begin{equation}\label{developtau}
\tau(\eta)-\bar \tau_1-2(\eta-\bar \eta_2)\ph_2(\gamma(\bar \eta_2))\,\det C=
2\,\det C\;\int_{\bar \eta_2}^\eta\,\Big(\ph_2\big(\gamma(s)\big)-\ph_2\big(\gamma(\bar \eta_2)\big)\Big)\;ds\,.
\end{equation}
Due to \cite{ASCV} and \cite{BSC2}, $\ph_2$ is a distributional solution
of \eqref{CcharDafappl}, hence Theorem~\ref{Dafermos} shows that
$\ph_2\big(s,\tau(s)\big)$ is continuously differentiable on $[a,b]$.
Then \eqref{developtau} leads us to the conclusion. $\Box$

\bigskip

\section{Selecting a flow of a continuous vector field}\label{selectflow}

In the next theorem we study the level set of a continuous function, assuming
a strict monotonicity on a family of curves that are solutions to a continuous
ODE, where there is no uniqueness.
\medskip
\begin{The}\label{teo:curva}
Let $A\subset\R^2$ be an open set with $(0,0)\in A$ and let $h:A\to \R$ be continuous.
Let $F:A\to \R$ be continuous, let $F(0,0)=0$ and assume that
each solution $\tau:I\lra\R$ of
\begin{equation}\label{ODE}
\taudot(\eta) = h(\eta,\tau(\eta))\,,
\end{equation}
whose graph is contained in $A$, has the property that 
$\eta\lra F(\eta,\tau(\eta))$ is strictly increasing on the
compact interval $I$.
Then there exists a compact neighbourhood of the origin $U\subset A$ and
an injective continuous curve $\zeta:[0,1]\to U$ such that
\begin{equation}\label{czeros}
\zeta([0,1])=U\cap F^{-1}(0).
\end{equation}
\end{The}
{\sc Proof.}
Let $R = [-a,a]\times [-b,b]$ be contained in $A$, set $M = \max_R |h|$ and 
$\delta = \min\{a, \frac{b}{2M}\}$.
If $|\tau(0)|\leq b/2$, then we have at least one
solution $\tau$ of \eqref{ODE} defined
on $I_\delta=[-\delta,\delta]$, whose graph is contained in $[-a,a]\times[-b,b]$.
Let $\bar \tau$ and $\hat \tau$ be, respectively, the infimum and the supremum of the solutions of \eqref{ODE} satisfying $\tau(0) = 0$.
Then $\bar \tau$ and $\hat \tau$ solve \eqref{ODE} and satisfy 
\[
\bar \tau \leq \hat \tau \quad\mbox{on}\quad I_\delta\quad
\mbox{and}\quad  \bar \tau(0) = \hat \tau(0) = 0. 
\]
Since $F(0,0) = 0$ and both $\eta\to F(\eta,\bar \tau(\eta))$,
$\eta\to F(\eta,\hat \tau(\eta))$ are strictly increasing,
whenever $\eta\in (0,\delta]$ we have $F(\eta,\bar \tau(\eta))>0$ and
$F(\eta,\hat \tau(\eta))>0$ and for every $\eta\in [-\delta,0)$ we have
$F(\eta,\bar \tau(\eta))<0$ and $F(\eta,\hat \tau(\eta))<0$.
We denote by $\tau_-$ and $\tau_+$ two solutions of \eqref{ODE} on $I_\delta$ satisfying 
\[
\tau_+(0) = +b/2,\quad\tau_-(0) = -b/2,\quad\mbox{and}\quad
\tau_-\leq \bar \tau\leq \hat \tau \leq \tau_+ \mbox{ on } I_\delta.
\]
Let us now apply Lemma~\ref{lemMajer} and then find two continuous and
nondecreasing one-parameter families of solutions to \eqref{ODE},
$\tau_\mu$ and $\tau_\nu$ defined on $[\hat\mu,\mu_+]$ and $[\nu_-,\bar\nu]$
that connect $\hat\tau$ to $\tau_+$ and $\tau_-$ to $\bar\tau$, where
\[
\hat\mu=\int_{I_\delta}\hat \tau,\quad
\mu_+=\int_{I_\delta} \tau_+,\quad 
\nu_-=\int_{I_\delta}\tau_-,\quad \bar\nu=\int_{I_\delta} \bar \tau.
\] 
By continuity of $[\hat\mu,\mu_+]\ni\mu\to\tau_\mu\in C(I_\delta,[-b,b])$
and $[\nu_-,\bar\nu]\ni\nu\to\tau_\nu\in C([I_\delta,[-b,b])$,
we can find $\nu_0<\bar\nu$ and $\mu_0>\hat\mu$
such that for $\xi\in[\nu_0,\bar\nu]\cup[\hat\mu,\mu_0]$, we have
\[
F(-\delta,\tau_\xi(-\delta))<0<F(\delta,\tau_\xi(\delta)).
\]
By strict monotonicity of $\eta\lra F(\eta,\tau_\xi(\eta))$,
for all $\xi\in[\nu_0,\bar\nu]\cup[\hat\mu,\mu_0]$,
it follows that there exists a unique point 
$\zeta=\zeta(\xi)\in\R^2$ on the graph of $\tau_\xi$ such that 
$F(\zeta(\xi)) = 0$.
Uniqueness of $\zeta(\xi)$ and the continuity of $\xi\to\tau_\xi
\in C(I_\delta,[-b,b])$ imply the continuity of
\[
\left[\nu_0,\bar\nu]\cup[\hat\mu,\mu_0\right]\ni\xi\lra\zeta(\xi)
\in I_\delta\times[-b,b]\,.
\]
Let us check that for all $\nu\in[\nu_-,\bar\nu[$, we have
$\tau_\nu(0)<\bar\tau(0)=0$.
The construction of $\tau_\nu$ in Lemma~\ref{lemMajer} gives
\[
\nu=\int_{I_\delta}\tau_\nu<\int_{I_\delta}\bar\tau=\bar\nu
\quad\mbox{and}\quad \tau_\nu\leq\bar\tau\,,
\]
hence $\tau_\nu\neq\bar\tau$.
From minimality of $\bar\tau$, it follows that
$\tau_\nu(0)<\bar\tau(0)=0$. One argues in a similar way to get 
$\tau_\mu(0)>\hat\tau(0)=0$ for all $\mu\in]\hat\mu,\mu_+]$.
This proves that $\tau_{\nu_0}(0)<0<\tau_{\mu_0}(0)$. 
Continuity of both $\tau_{\nu_0}$ and $\tau_{\mu_0}$ implies that the compact set 
\[
U=\{(\eta,\tau)\in\R^2:\ \eta\in I_\delta,\, \tau_{\nu_0}(\eta)\leq\tau\leq\tau_{\mu_0}(\eta)\}
\]
is a neighbourhood of $(0,0)$.
Take any point $(\eta_0,\tau_0)\in U$ and assume that
\[
\bar\tau(\eta_0)\leq\tau_0\leq\hat\tau(\eta_0).
\]
If $\eta_0=0$, then $(\eta_0,\tau_0)=(0,0)$ is a zero of $F$.
If $\eta_0>0$, we can choose a maximal solution $y$ of \eqref{ODE}
on some open interval $J=]\alpha,\beta[$ such that $y(\eta_0)=\tau_0$.
By maximality of $y$, if $\beta\leq\delta$, then there exists $\ep>0$ such that
$y(t)<\bar\tau(t)$ for all $t\in[\beta-\ep,\beta[$.
Then $\max\{y(t),\bar\tau(t)\}=\bar\tau(t)$ for all
$t\in[\beta-\ep,\beta[$, therefore this function is a solution to
\eqref{ODE} that can be extended to 
$]\max\{\alpha,-\delta\},\delta]$.
We argue in a similar way in the case $-\delta\leq\alpha$, hence 
$\max\{y(t),\bar\tau(t)\}$ is well defined on $I_\delta$
and satisfies \eqref{ODE}.
Now, we set $\tilde\tau=\min\{\hat\tau,\max\{y,\bar\tau\}\}$.
We observe that $\tilde\tau$ is still a solution of \eqref{ODE} and 
clearly $\tilde\tau(\eta_0)=\tau_0$ and
$\bar\tau\leq\tilde\tau\leq\hat\tau$.
Since $\tilde\tau(0)=0$ and $\eta\to F(\eta,\tilde\tau(\eta))$ 
is strictly increasing, we have
$F(\eta_0,\tilde\tau(\eta_0))=F(\eta_0,\tau_0)>0$.
If $\eta_0<0$, then one argues in the same way getting $F(\eta_0,\tau_0)<0$.

By construction of $\zeta$, we observe that $\zeta:\left[\nu_0,\bar\nu]\cup[\hat\mu,\mu_0\right]\lra U$ and also
\[
\zeta\big(\left[\nu_0,\bar\nu]\cup[\hat\mu,\mu_0\right]\big)
\subset U\cap F^{-1}(0)\,.
\]   
Let us now pick $(\eta_0,\tau_0)\in U\cap F^{-1}(0)$. 
If $(\eta_0,\tau_0)=(0,0)$, then it coincides with both
$\zeta(\bar\nu)$ and $\zeta(\hat\mu)$, so that we can assume
$(\eta_0,\tau_0)\neq(0,0)$.
The previous arguments imply that either 
$\tau_0<\bar\tau(\eta_0)$ or $\tau_0>\hat\tau(\eta_0)$.
Let us consider for instance the case $\tau_0>\hat\tau(\eta_0)$.
Since $(\eta_0,\tau_0)\in U$, we have 
\[
\tau_{\hat\mu}(\eta_0)=\hat\tau(\eta_0)<\tau_0\leq\tau_{\mu_0}(\eta_0),
\]
hence the number
\[
\mu_1=\inf\{\mu\in[\hat\mu,\mu_0]: \tau_\mu(\eta_0)\geq\tau_0\}
\]
is well defined and by continuity of $\mu\to\tau_\mu(\eta_0)$, we get
\[
\hat\mu<\mu_1\leq\mu_0\quad\mbox{and}\quad
\tau_{\mu_1}(\eta_0)\geq\tau_0.
\]
By contradiction, if we had $\tau_{\mu_1}(\eta_0)>\tau_0$,
again continuity of $\mu\to\tau_\mu(\eta_0)$ would contradict the definition of $\mu_1$ giving $\mu_1'\in]\hat\mu,\mu_1[$ such that 
$\tau_{\mu_1'}(\eta_0)>\tau_0$. This shows that
$\tau_{\mu_1}(\eta_0)=\tau_0$, hence we have $\eta_0\in I_\delta$ such that 
$F(\eta_0,\tau_{\mu_1}(\eta_0))=0$, by the strict monotonicity of
$\eta\to F(\eta,\tau_{\mu_1}(\eta))$ we must have 
\[
\zeta(\mu_1)=(\eta_0,\tau_{\mu_1}(\eta_0))=(\eta_0,\tau_0)\,.
\]
We have proved that $(\eta_0,\tau_0)\in
\zeta\big(\left[\nu_0,\bar\nu]\cup[\hat\mu,\mu_0\right]\big)$.
Clearly, in the case $\tau_0<\bar\tau(\eta_0)$ the same previous
argument holds, hence leading us to the equality
\[
\zeta\big(\left[\nu_0,\bar\nu]\cup[\hat\mu,\mu_0\right]\big)
=U\cap F^{-1}(0)\,.
\]   
Since $\zeta(\bar\nu)=\zeta(\hat\mu)=(0,0)$, an elementary change
of variable allows us to assume that $\zeta:[0,1]\lra U\cap F^{-1}(0)$.
is a surjective continuous curve.
To conclude, we have to show that $\zeta$ can be reparameterized as an injective, continuous curve. We just notice that 
$\zeta(\xi_1) = \zeta(\xi_2)$ implies $\zeta(\xi)= \zeta(\xi_1)$ for all $\xi_1\leq \xi\leq \xi_2$ due to monotonicity of $\xi\to\tau_\xi$ in the sense of 
Lemma \ref{lemMajer}. This property of $\zeta$ corresponds to the fact
that all of its preimages are intervals. It is an elementary fact to notice
that this property gives the existence of an injective
continuous reparametrization of $\zeta$. $\Box$

\medskip
\noindent
The proof of the previous theorem requires the following lemma,
that we have not found in the existing literature on
ODEs and that has been kindly pointed out to us by P. Majer.
\medskip
\begin{Lem}\label{lemMajer}
Let $A$ be an open set of $\R^2$,
let $I$ and $J$ be two compact intervals such that $I\times J\subset A$
and let $\tau_-$, $\tau_+$ be two solutions of \eqref{ODE} defined on $I$,
whose graphs are contained in $I\times J$ and such that 
$\tau_-\leq \tau_+$ on $I$. Then defining the numbers
\[
\mu_- = \int_I \tau_-(\eta)\,d\eta\quad\mbox{and}\quad
\mu_+ = \int_I \tau_+(\eta)\, d\eta,
\]
there exists a continuous curve $[\mu_-,\mu_+]\ni\mu\to \tau_\mu$
with respect to the $L^\infty$-norm on the space of solutions to \eqref{ODE}
defined on $I$, such that $\tau_{\mu_-} = \tau_-$ and $\tau_{\mu_+} = \tau_+$.
Moreover, $\tau$ satisfies
\begin{itemize}
\item[(1)] if $\mu_1\leq \mu_2$, then $\tau_{\mu_1}\leq\tau_{\mu_2}$ on $I$,
\item[(2)] $\displaystyle\int_I \tau_\mu(\eta)\, d\eta = \mu$ for all $\mu\in [\mu_-,\mu_+]$.
\end{itemize}
\end{Lem}
{\sc Proof.}
Set $R = I\times J$ and let $S$ be the family of solutions to \eqref{ODE},
defined on $I$ such that $\tau_-\leq \tau \leq \tau_+$. 
Clearly, the graph of solutions of $S$ is contained in $R$.
For any couple of solutions $\tau,\rho\in S$ such that $\rho\leq\tau$ on $I$
we introduce the set 
\[ 
\ccI(\rho,\tau) = \{\tau\in S:\  \rho\leq \tau\leq\tau\; \mbox{on $I$}\}\,.
\]
We first notice that for every $\tau_1,\tau_2\in S$ such that $\tau_1\leq \tau_2$ on $I$ the 
family $\ccI(\tau_1,\tau_2)$ is connected. 
Taking into account that each $\tau\in S$ has a 
Lipschitz constant less than or equal to $M = \max_{R} |h(\eta,\tau)|$,
then the proof of connectedness follows the same arguments used to prove connectedness
of the interval. For the sake of the reader, we sketch here a few details.
By contradiction, let $C_1$ and $C_2$ be disjoint compact sets of $S$
such that $\ccI(\tau_1,\tau_2)\subset C_1\cup C_2$ and $C_1$
does not contain $\tau_2$. 
Let $\rho_1(\eta)=\sup_{\tau\in C_1} \tau(\eta)$ be the upper envelope
of $C_1$.  Taking into account that both pointwise maximum of two solutions
is still a solution, one can find a converging sequence of solutions in $C_1$
that restricted on a dense subset of $I$ pointwise converge to $\rho_1$.
This shows that $\rho_1\in C_1$. Similarly, we set
\[
\rho_2(\eta)=\inf_{\tau\in C_2\cap\ccI(\rho_1,\tau_2)}\tau(\eta)
\]
and observe that it belongs to $C_2\cap\ccI(\rho_1,\tau_2)$.
Clearly, the lower envelope $\rho_2\geq\rho_1$ cannot coincide with $\rho_1$,
since $\rho_1\in C_1$.
Then there exists $a\in I$ where $\rho_1(a)<\rho_2(a)$ and one
can find a solution $\tau_0$ of $S$ that satisfies 
$\rho_1(a)<\tau_0(a)<\rho_2(a)$.
We define 
\[
\bar\tau=\max\big\{\rho_1,\min\{\tau_0,\rho_2\}\big\}
\]
that is clearly still a solution of (\ref{ODE}),
satisfies $\rho_1\leq\bar\tau\leq\rho_2$ and $\bar\tau(a)=\tau_0(a)$.
In particular, $\bar\tau\in\ccI(\tau_1,\tau_2)$ and it 
differs from both $\rho_1$ and $\rho_2$.
This leads us to a contradiction with the inclusion
$\ccI(\tau_1,\tau_2)\subset C_1\cup C_2.$
In view of connectedness of $\ccI(\tau_1,\tau_2)$, we have  
\begin{equation}\label{muconnect}
 \left\{\int_I \tau:\ \tau\in {\mathfrak I}(\tau_1,\tau_2)\right\} 
= \left[\int_I \tau_1, \int_I \tau_2\right].
\end{equation}
The previous property, through a dyadic construction allows us
to obtain a continuous section of the function
$\cL:\ccI(\tau_-,\tau_+)\lra[\mu_-,\mu_+]$, defined as
$\cL(\tau)=\int_I\tau(\eta)\,d\eta$.
We first choose $\mu_{1,1} = \frac{\mu_- + \mu_+}{2}$, then thanks to 
\eqref{muconnect} we select
 $\tau_{\mu_{1,1}} \in {\mathfrak I}(\tau_-,\tau_+)$ such that
\[
 \int_I \tau_{\mu_{1,1}} = \mu_{1,1}.
\]
We then proceed on the two subintervals $[\mu_-,\mu_{1,1}]$ and $[\mu_{1,1},\mu_+]$, selecting their middle points $\mu_{2,1}$ and $\mu_{2,2}$, respectively, and choosing $\tau_{\mu_{2,1}}$ and $\tau_{\mu_{2,2}}$ satisfying
\[
 \tau_- \leq \tau_{\mu_{2,1}} \leq \tau_{\mu_{1,1}}, \qquad 
 \tau_{\mu_{1,1}} \leq \tau_{\mu_{2,2}} \leq \tau_+\quad
\mbox{and}\quad
 \int_I \tau_{\mu_{2,i}} = \mu_{2,i}\quad\mbox{for}\quad i=1,2.
\]
Iterating this procedure, we get a nondecreasing mapping $\nu\to \tau_\nu$ 
defined on a dense subset of $[\mu_-,\mu_+]$, taking values in $\ccI(\tau_-,\tau_+)$.
For every $\mu\in[\mu_-,\mu_+]$, we define the extension
\[
\tau_\mu(\eta)=\lim_{\stackrel{\nu\to\mu^-}{\nu\mbox{\tiny\ dyadic}}}\tau_\nu(\eta)
\]
for all $\eta\in I$. Compactness of $\ccI(\tau_-,\tau_+)$ makes
pointwise converging sequences into uniform converging sequences,
up to subsequences. Then $[\mu_-,\mu_+]\ni\mu\lra\tau_\mu\in\ccI(\tau_-,\tau_+)$
is clearly a continuonus section of $\cL$ and satisfies the claimed properties. $\Box$

\section{Proof of Theorem~\ref{ImplCurv} and remarks on regularity} 

Collecting all results of the previous sections, we are now able to prove
Theorem~\ref{ImplCurv}. First of all, up to left translations, it is not restrictive
to assume that $0\in\Sigma_1\cap\Sigma_2$. The $\H$-regular surfaces
$\Sigma_1$, $\Sigma_2$ are represented in a suitable neighbourhood
$\cU$ of $0$ as the zero level sets of $f_1,f_2\in\cC^1(\cU,\R)$, respectively.
In addition, we can assume that $\nabla_Hf_1$ and $\nabla_Hf_2$ are 
linearly independent on $\cU$, since the horizontal normals of $\Sigma_1$ and 
$\Sigma_2$ at the origin are linearly independent.
We set $f=(f_1,f_2)$ and observe that
\[
\Sigma_1\cap\Sigma_2\cap\cU=f^{-1}(0)\,.
\]
Since $\nabla_Hf(x)$ is surjective for all $x\in\cU$, it is not restrictive
to assume that $Y_1f_2(0)\neq0$.
We will follow the assumptions of Section~\ref{DirectDerivChRu},
setting $Y_1=c^1_1X_1+c^2_1X_2$ and $b_1=c^1_1e_1+c^2_1e_2$.
Then $N$ is a vertical subgroup of orthonormal basis $(b_2,e_3)$, 
such that $H\oplus N=\H$, where $H=\span\{b_1\}$.
The implicit function theorem of \cite{FSSC3} gives us an open set
$U$ of $N$ along with a continuous mapping $\phi_2:U\lra H$ such that,
up to further shrinking $\cU$, we have $\Phi_2(U)=\Sigma_2\cap\cU$,
where $\Phi_2(n)=n\cdot\phi_2(n)$ and $\phi_2(n)=\ph_2(n)b_1$.
Then the intersection can be written as follows
\[
\Sigma_1\cap\Sigma_2\cap\cU=\{\Phi_2(n)\in\cU: f_1\big(\Phi_2(n)\big)=0\}\,.
\]
We identify $U\subset N$ with an open set $A\subset\R^2$ through
the basis $(b_2,e_3)$ and consider $F=f_1\circ\Phi_2:A\lra\R$,
observing that $F(0,0)=0$. Let us now consider any characteristic $\tau$ 
of \eqref{CcharDafappl}, namely, a solution of 
\[
\dot\tau(\eta)=2\det(C)\,\ph_2\big(\eta,\tau(\eta)\big)\,.
\]
Then we can apply Theorem~\ref{chainr}, getting
that $F\big(\eta,\tau(\eta)\big)$ is differentiable and
\[
\frac{d}{d\eta}F\big(\eta,\tau(\eta)\big)=
-\frac{1}{Y_1f_2\big(\Phi_2\big(\eta,\tau(\eta)\big)\big)}\;\det
\left(\begin{array}{cc}
Y_1f_1\big(\Phi_2\big(\eta,\tau(\eta)\big)\big) & Y_2f_1\big(\Phi_2\big(\eta,\tau(\eta)\big)\big) \\ 
Y_1f_2\big(\Phi_2\big(\eta,\tau(\eta)\big)\big) & Y_2f_2\big(\Phi_2\big(\eta,\tau(\eta)\big)\big)
\end{array}\right)\neq0\,,
\]
whenever $\big(\eta,\tau(\eta)\big)\in A$. In particular, 
$\eta\to F\big(\eta,\tau(\eta)\big)$ is strictly monotone.
This property allows us to apply Theorem~\ref{teo:curva}
with $h=2\det(C)\,\ph_2$, hence getting a neighbourhood
$U\subset A$ of the origin such that 
$U\cap F^{-1}(0)$ is the image of an injective continuous
curve $\zeta:[0,1]\lra U$. 
Then we can find an open set $O\subset\cU$ such that 
\[
\Sigma_1\cap\Sigma_2\cap O=\{\Phi_2(\zeta(\xi))\in O: \xi\in[0,1]\}\,.
\]  
Setting $\Gamma=\Phi_2\circ\zeta$, the proof of Theorem~\ref{ImplCurv} is concluded. $\Box$  

\subsection{A cone-type Lipschitz continuity}\label{intLip}

In analogy with intrinsic cones of \cite{FSSC7}, associated to a 
semidirect factorization of $\H^n$, we introduce similar cones in $\H$,
although here there is no semidirect factorization.
Recall the canonical decomposition associated to any element
$x=x_1+x_2\in\H$, where $x_j\in H_j$ and $j=1,2$. 
The difference with respect to a semidirect factorization is that 
here $H_1$ is not a subgroup.

Given $a,r>0$ and $p\in \H^n$, we define the (intrinsic)
{\em closed cone} with base $H_2$, axis $H_1$, width $r>0$
and opening $\alpha>0$ as
\[
C_r(\alpha) =\big\{x\in\H: \|x_2\|\leq \alpha\, \|x_1\|\leq\alpha\,r \big\}.
\]
The closed cone with vertex $p$ is the translated cone
\[
C_r(p,\alpha)=p\cdot C_r(\alpha)\,.
\]
A set $S\subset\H$ has the {\em cone property} if for every $p\in S$
we can find a neighbourhood $U$ of $p$ such that for all 
$\alpha>0$ there exist $r>0$, depending on $\alpha$ and $U$,
such that for all $x\in U\cap S$ there holds
\begin{equation}
\label{Hcond}
S\cap C_r(x,\alpha) = \{x\} \,.
\end{equation}
It is easy to observe that any level set of
$f\in\cC^1(\Omega,\R^2)$, where $\nabla_Hf$ is everywhere
surjective has the cone property. 
In fact, let $p\in f^{-1}(z)$ and let $\delta>0$ and $\omega_p$ be as in
Lemma~\ref{unidiff}. By surjectivity, we can make $\delta$
small such that
\[
\lambda=\min_{\substack{x'\in D_{p,\delta}\\ v\in H_1,\,\|v\|=1}}
|\nabla_Hf(x')(v)|>0\,.
\]
Then for any $x,y\in D_{p,\delta/2}\cap f^{-1}(z)$ we have
\begin{equation}\label{angleest}
\lambda\,|y_1-x_1|
\leq|\nabla_Hf(x)(x^{-1}\cdot y)|\leq \omega_p(\|x^{-1}\cdot y\|)\,
\|x^{-1}\cdot y\|\,.
\end{equation}
Let $\alpha>0$ and set $0<\ep<\lambda/(\alpha+1)$. Let $t_\ep>0$ be
such that $\displaystyle{ \sup_{0\leq s\leq t_\ep}\omega_p(s)<\ep}$
and $t_\ep<\delta/4$. Thus, for every $x\in D_{p,\delta/4}\cap f^{-1}(z)$ 
and $y\in D_{x,t_\ep}\cap f^{-1}(z)\sm\{x\}$, 
as a consequence of \eqref{angleest} we get
\[
|y_1-x_1|<\frac{1}{\alpha}\,\sqrt{|(x^{-1}\cdot y)_2|}\,,
\]
where $x^{-1}\cdot y=(x^{-1}\cdot y)_1+(x^{-1}\cdot y)_2$ and $(x^{-1}\cdot y)_j\in H_j$.
Finally, defining $r=t_\ep/\max\{1,\alpha\}$, we have proved that 
$C_r(x,\alpha)\cap f^{-1}(z)=\{x\}$.


\begin{thebibliography}{99}
%
%
%
\bibitem{AmbKir00}
{\sc L.Ambrosio and B.Kirchheim},
{\em Rectifiable sets in metric and {B}anach spaces}, 
Math. Ann., {\bf 318}, n.3, 527-555, (2000).

\bibitem{ASCV}{\sc L.Ambrosio, F.Serra Cassano, D.Vittone},
{\em Intrinsic regular hypersurfaces in Heisenberg groups},
J. Geom. Anal. {\bf 16}, no. 2, 187-232, (2006)

\bibitem{AreSer09}{\sc G.Arena, R.Serapioni},
{\em Intrinsic regular submanifolds in Heisenberg groups 
are differentiable graphs}, Calc. Var., {\bf 35}, 517-536, (2009)

\bibitem{Bell}{\sc A.Bella\"iche},
{\em The Tangent space in sub-Riemannian geometry}, in {\em Subriemannian
Geometry}, Progress in Mathematics, {\bf 144}. ed. by A.Bellaiche and
J.Risler, Birkh\"auser Verlag, Basel, (1996).

\bibitem{BSC1}{\sc F.Bigolin, F.Serra Cassano},
{\em Intrinsic regular graphs in Heisenbegr groups vs. 
weak solutions of non linear first-order PDEs}, to appear in
Adv. Calc. Var

\bibitem{BSC2}{\sc F.Bigolin, F.Serra Cassano},
{\em Distributional solutions of Burgers' equation and intrinsic
regular graphs in Heisenberg groups}, to appear in J. Math. Anal. Appl.

\bibitem{Dafermos06}{\sc C.M.Dafermos},
{\em Continuous solutions for balance laws},
Ric. Mat. {\bf 55}, n.1, 79-91, (2006)

\bibitem{FS82}{\sc G.B.Folland, E.M. Stein},
{\em Hardy Spaces on Homogeneous groups}, Princeton University Press, (1982)

\bibitem{FSSC3}{\sc B.Franchi, R.Serapioni, F.Serra Cassano},
{\em Rectifiability and Perimeter in the Heisenberg group},
Math. Ann. {\bf 321}, n.3, 479-531, (2001)

\bibitem{FSSC6}{\sc B.Franchi, R.Serapioni, F.Serra Cassano},
{\em Regular submanifolds, graphs and area formula in Heisenberg
groups}, Adv. Math. {\bf 211}, n.1, 152-203, (2007)

\bibitem{FSSC7}{\sc B.Franchi, R.Serapioni, F.Serra Cassano},
{\em Intrinsic Lipschitz graphs in Heisenberg groups}, 
J. Nonlinear Convex Anal. , {\bf 7}, n.3, p. 423-441, (2006)

\bibitem{KirSC04}
{\sc B.Kirchheim and F.Serra Cassano},
{\em Rectifiability and parameterization of intrinsic regular surfaces in
 the Heisenberg group},
Ann. Sc. Norm. Super. Pisa Cl. Sci. (5), 3(4), 871-896 (2004).

\bibitem{Mag15}{\sc V.Magnani},
{\em Area implies coarea}, Indiana Univ. Math. J., to appear, (2010)  

\bibitem{Rum90}{\sc M.Rumin},
{\em Un complexe de formes diff\'erentielles sur les vari\'et\'es de contact,}
C. R. Acad. Sci. Paris S\'er. I Math. {\bf 310}, n.6, 401-404, (1990)

%
%
%
\end{thebibliography}
\end{document}